\documentclass[11pt]{amsart}
\usepackage{amsxtra}
\usepackage{amssymb}
\theoremstyle{plain}
\newcommand{\cleqn}{\setcounter{equation}{0}}
\newcommand{\clth}{\setcounter{theorem}{0}}
\newcommand {\sectionnew}[1]{\section{#1}\cleqn\clth}
\newcommand{\nn}{\hfill\nonumber}
\newtheorem{theorem}{Theorem}[section]
\newtheorem{lemma}[theorem]{Lemma}
\newtheorem{definition-theorem}[theorem]{Definition-Theorem}
\newtheorem{proposition}[theorem]{Proposition}
\newtheorem{corollary}[theorem]{Corollary}
\newtheorem{conjecture}[theorem]{Conjecture}
\newtheorem{definition}[theorem]{Definition}
\newtheorem{example}[theorem]{Example}
\newtheorem{remark}[theorem]{Remark}
\newtheorem{notation}[theorem]{Notation}
\newcommand \bth[1] { \begin{theorem}\label{t#1} }
\newcommand \ble[1] { \begin{lemma}\label{l#1} }

\newcommand \bpr[1] { \begin{proposition}\label{p#1} }
\newcommand \bco[1] { \begin{corollary}\label{c#1} }
\newcommand \bde[1] { \begin{definition}\label{d#1}\rm }
\newcommand \bex[1] { \begin{example}\label{e#1}\rm }
\newcommand \bre[1] { \begin{remark}\label{r#1}\rm }
\newcommand \bcon[1] { \begin{conjecture}\label{cn#1}\rm}

\newcommand \bnota[1] { \begin{notation}\label{n#1}\rm }
\newcommand {\ethe} { \end{theorem} }
\newcommand {\ele} { \end{lemma} }

\newcommand {\epr} { \end{proposition} }
\newcommand {\eco} { \end{corollary} }
\newcommand {\econ} { \end{conjecture} }
\newcommand {\ede} { \end{definition} }
\newcommand {\eex} { \end{example} }
\newcommand {\ere} { \end{remark} }

\newcommand {\enota} { \end{notation} }
\newcommand \thref[1]{Theorem \ref{t#1}}
\newcommand \leref[1]{Lemma \ref{l#1}}
\newcommand \prref[1]{Proposition \ref{p#1}}
\newcommand \coref[1]{Corollary \ref{c#1}}

\newcommand \reref[1]{Remark \ref{r#1}}
\newcommand \lb[1]{\label{#1}}
\def \Rset {{\mathbb R}}         %mathsets
\def \Cset {{\mathbb C}}
\def \Zset {{\mathbb Z}}

\def \O  {{\mathcal{O}}}

\def \al {\alpha}

\def \la {\lambda}

\def \ol {\overline}
\def \wt {\widetilde}

\def \Id { {\mathrm{Id}} }

\def \Span { {\mathrm{Span}} }

\def \rank { {\mathrm{rank}} }

\def \gl  {\mathfrak{gl}}

\def \Y {\Upsilon}
\def \dx {\partial_x}
\DeclareMathOperator \ad { {\mathrm{ad}} }

\DeclareMathOperator \Gr { {\mathrm{Gr}} }
\DeclareMathOperator \Hom { {\mathrm{Hom}} }

\DeclareMathOperator \End { {\mathrm{End}} }
\DeclareMathOperator \Spec { {\mathrm{Spec}} }
\DeclareMathOperator \diag { {\mathrm{diag}} }

\DeclareMathOperator \Wr { {\mathrm{Wr}} }
\DeclareMathOperator \tr { {\mathrm{tr}} }
\DeclareMathOperator \reg { \sf reg }
\renewcommand \Im { {\mathrm{Im}} \, }
\renewcommand \Re { {\mathrm{Re}} \, }

\newcommand{\Rr}{\Rset^n_{\sf reg}}
\newcommand{\Cr}{\Cset^n_{\sf reg}}
\newcommand{\RrCS}{\Rset^{(n)}_{\sf reg}}
\newcommand{\CrCS}{\Cset^{(n)}_{\sf reg}}
\newcommand{\Mat}{{\mathrm{Mat}}}

\begin{document}
%%%%%%%%%%%%%%%%%%%%%%    Title    %%%%%%%%%%%%%%%%%%%%%%%%%%%%%%%%%%%%%%%%
\title[Calogero--Moser spaces and Cherednik algebras]
{The real loci of Calogero--Moser spaces, representations
of rational Cherednik algebras and the Shapiro conjecture}
\author[Iain Gordon]{Iain Gordon}
\address{I.G.: School of Mathematics and Maxwell Institute for 
Mathematical Sciences, Edinburgh University, Edinburgh EH9 3JZ, Scotland}
\email{igordon@ed.ac.uk}
\author[Emil Horozov ]{Emil Horozov}
\address{E.H.: Department of Mathematics and Informatics,
Sofia University, 5 J. Bourchier Blvd., Sofia 1126, Bulgaria, and }
\address{Institute of Mathematics and Informatics, Bulg. Acad. of Sci.,
Acad. G. Bonchev Str., Block 8, 1113 Sofia, Bulgaria }
\email{horozov@fmi.uni-sofia.bg}
\author[Milen Yakimov]{Milen Yakimov}
\address{M.Y.: Department of Mathematics \\
University of California \\
Santa Barbara, CA 93106, U.S.A.} \email{yakimov@math.ucsb.edu}
%\thanks{}
\date{}
%\keywords{Calogero--Moser spaces, the Shapiro--Shapiro conjecture, 
%Wilson's adelic Grassmannian, 
%Cherednik algebras.}
%\subjclass[2000]{Primary 53D17; Secondary 58H05, 17B62}
\begin{abstract}
We prove a criterion for the reality of irreducible representations
of the rational Cherednik algebras $H_{0,1}(S_n)$. This is shown
to imply a criterion for the real loci
of the Calogero--Moser spaces $C_n$ in terms of the 
Etingof--Ginzburg finite maps 
$\Upsilon \colon C_n \to \Cset^n/S_n \times \Cset^n/S_n$,
recovering a result 
of Mikhin, Tarasov, and Varchenko \cite{MTV2}.
As a consequence we obtain a criterion 
for the real locus of the Wilson's adelic Grassmannian of rank one bispectral 
solutions of the KP hierarchy. 
Using Wilson's first parametrisation of the adelic 
Grassmannian, we give a new proof of a result of \cite{MTV2}
on real bases of spaces of quasi polynomials. The Shapiro 
Conjecture for Grassmannians is equivalent to a special case
of our result for Calogero--Moser spaces, namely for 
the fibres of $\Upsilon$ over $\Cset^n/S_n \times 0$.     
\end{abstract}
\maketitle
%%%%%%%%%%%%%%%%%%%%   Introduction   %%%%%%%%%%%%%%%%%%%%%%%%%%%%%%%%%%%%%%%%
\sectionnew{Introduction}\lb{intro}
The $n$-th Calogero--Moser space $C_n$ is the geometric quotient 
of 
\[  
\ol{C}_n = \{ (X, Z) \in \gl_n(\Cset)^{\times 2} \mid 
\rank ([X,Z] + I_n)=1 \}
\]
by the action of $GL_n(\Cset)$ by simultaneous conjugation. 
It is a smooth, irreducible, complex, affine variety, \cite{W2}. 
The space $C_n$ is the phase space of the (complex) 
Calogero--Moser integrable system \cite{KKS,W2}
and parametrizes irreducible representations 
of the deformed preprojective algebra of a certain quiver 
\cite{CBH,CB}. We define the real locus $RC_n$ of $C_n$ as the 
image under $\pi_n$ of the space of pairs of real matrices 
inside $C_n$. It is not hard to see that $RC_n$ is a real algebraic subset
of $C_n$ which is isomorphic to the $n$-th real Calogero--Moser 
space.
 
A different interpretation of the spaces $C_n$ 
was found by Etingof and Ginzburg \cite{EG} in terms of representations
of rational Cherednik algebras associated to symmetric groups. 
%Among many other applications, they succeeded to construct 
%Calogero--Moser spaces for all Coxeter groups. 

In this paper we show how to use the representation theory of rational 
Cherednik algebras to obtain results on the real algebraic geometry
of $C_n$. In particular we give new proofs of several theorems of Mukhin, Tarasov and Varchenko, \cite{MTV1, MTV2}, including the Shapiro Conjecture for Grassmannians. 

\medskip
The rational Cherednik algebra $H_{0,1}(S_n)$ is a specialisation of 
a two parameter deformation of the smash product 
$ \Cset[x_1, \ldots, x_n, y_1, \ldots, y_n]\rtimes \Cset S_n$.
Its irreducible representations all have complex dimension
$n!$ and are parametrised by the points of $C_n$, \cite{EG}.
Denote by $\ol{e}$ the symmetrising idempotent of the 
copy of $S_{n-1}$ inside $S_n$ which permutes only the last 
$n-1$ indices. Given an irreducible $H_{0,1}(S_n)$ module 
$V$, $x_1$ and $y_1$ preserve the $n$-dimensional 
subspace $\ol{e} V$ and define a point 
$\pi_n(x_1|_{\ol{e}V},y_1|_{\ol{e}V}) \in C_n$.
Etingof and Ginzburg proved that this establishes a bijection 
between the equivalence classes of irreducible representations of 
$H_{0,1}(S_n)$ and the points of $C_n$. In order to state the main result 
of this paper, we note that $H_{0,1}(S_n)$ has a natural real form:
the real subalgebra $H_{0,1}^{\Rset}(S_n)$ generated by 
the elements of $S_n$ and $x_1, \ldots, x_n, y_1, \ldots, y_n$.

\medskip
\noindent {\bf Main Theorem} (\thref{Cher}). {\it If an irreducible representation $V$ of $H_{0,1}$
has the property that $x_1|_{\ol{e}V}$ and $y_1|_{\ol{e}V}$
have only real eigenvalues, then $V$ is the complexification 
of a real representation of $H_{0,1}^{\Rset}(S_n)$.}

\medskip
Now define the real locus $RC_n$ of $C_n$ as the image of the subset 
of real matrices 
in $\ol{C}_n$ -- this is nothing but the real locus of $\ol{C}_n$ --
under the quotient map $\pi_n \colon \ol{C}_n \to C_n$. We show that 
$RC_n$ is a real algebraic subset of $C_n$ which is isomorphic 
as a real affine variety to the $n$-th real Calogero-Moser space.

Set $\Cset^{(n)}= \Cset^n/S_n$ and let $\Spec(X)$ stand for the
eigenvalues of a square matrix $X$. In \cite{EG} Etingof and Ginzburg proved that the canonical map
\begin{equation}
\Y \colon C_n \to \Cset^{(n)} \times \Cset^{(n)}, \quad 
\Y (\pi_n(X,Z))= (\Spec(X), \Spec(Z)),
\label{map_Ups}
\end{equation}
is a finite map of degree $n!$. This map, and particularly its fibre over $0\times 0$, was studied in \cite{EG, FG, G}. We obtain from the Main Theorem

\medskip
\noindent {\bf Corollary} (\thref{2}). {\it Let $\Rset^{(n)}:= \Rset^n/S_n \subset \Cset^{(n)}$. We have \begin{equation}
\Y^{-1}( \Rset^{(n)} \times \Rset^{(n)}) 
\subset RC_n.
\label{Up1}
\end{equation}
}

In elementary terms this claims that
if $(X, Z) \in \ol{C}_n$ and both $X$ and $Z$ have real 
eigenvalues, then $X$ and $Z$ can be 
simultaneously conjugated (under $GL_n(\Cset)$) to 
pair of real matrices. This reproves a result of Mukhin, Tarasov and
Varchenko.

\medskip
The Calogero--Moser space $C_n$ parametrises the equivalence classes 
of representations of a specific dimension vector of the deformed preprojective algebra 
$\Pi_{\nu}(Q)$ of a certain quiver, \cite{CBH, CB}. As an immediate consequence of the Main Theorem,
we also obtain a criterion for reality of the representations of
$\Pi_{\nu}(Q)$ in this class.

\medskip
The disjoint union of all Calogero--Moser spaces also parametrizes
Wilson's adelic Grassmannian $\Gr^{\ad}$, \cite{W2}. The latter space first 
arose as the set of all solutions of the KP hierarchy which have 
bispectral wave functions of rank 1, \cite{W1}. We define and 
study in detail the real locus of $\Gr^{\ad}$. All possible 
approaches to the definition of the real locus of $\Gr^{\ad}$
(as the union of the real loci of $C_n$, or by requiring 
reality of the associated tau or wave functions) are shown to 
be equivalent. From the Main Theorem we derive the following
criterion:

\medskip
\noindent {\bf Corollary} (\thref{Wad}). {\it If $W \in \Gr^{\ad}$ has the property that the specialisations 
of the tau function $\tau_W(x, 0, \ldots)$ and the bispectral 
dual tau function $\tau_{b W}(x, 0, \ldots)$ have real roots, 
then $\tau_W(t_1, t_2, \ldots) \in \Cset[[t_1, t_2, \ldots]]$
has real coefficients.
}

\medskip
Wilson first defined the  
adelic Grassmannian (which actually motivated the term) 
by imposing a set of linear conditions of a special type 
on the plane of the 
trivial solution of the KP hierarchy. Translating 
the criterion for the real locus of $\Gr^{\ad}$ in terms of 
these conditions leads to another proof of the following result 
of Mukhin, Tarasov, and Varchenko. 

\medskip
\noindent {\bf Corollary} (\thref{3}). {\it Fix a collection of 
distinct real numbers $\mu_1, \ldots, \mu_k \in \Rset$ 
and a collection of  finite dimensional 
subspaces $V_1, \ldots, V_k$ of $\Cset[x]$. If for a given basis
$\{ q_1(x), \ldots, q_N(x) \}$ 
of $e^{\mu_1 x} V_1 \oplus \cdots \oplus e^{\mu_k x} V_k$
the polynomial
\[
e^{-(\mu_1 + \cdots + \mu_n) x} \Wr(q_1(x), \ldots, q_N(x))
\]
has only real roots, then all vector spaces $V_1, \ldots, V_k$
have real bases.}

\medskip
This paper was motivated by an attempt to understand 
the relationship between the Shapiro Conjecture, the 
Calogero--Moser spaces and the rational Cherednik algebra 
$H_{0,1}(S_n)$. The Shapiro Conjecture for Grassmannians is a consequence of the Main Theorem.

\medskip
\noindent {\bf Corollary.} {\it If $p_1(x), \ldots p_n(x) \in \Cset[x]$ are such that the Wronskian
$\Wr(p_1(x), \ldots, p_n(x))$ has only real roots, then 
$\Span \{ p_1(x), \ldots,$ $p_n(x) \}$ has a real basis.}

\medskip 
This result plays a major role in the real Schubert calculus \cite{S1,S2} and 
the theory of real algebraic curves \cite{KS}. Considerable
numerical evidence to support the conjecture was first obtained in \cite{S3}.
The conjecture was proved in the case $n=2$ by Eremenko and Gabrielov 
\cite{EG1, EG2}. In the general case, it was 
proved by Mukhin, Tarasov, and Varchenko \cite{MTV1}, who also proved 
the above generalisation of the conjecture for quasipolynomials \cite{MTV2}.
The Shapiro Conjecture is in fact equivalent to the special 
case of \eqref{Up1}
$
\Upsilon^{-1}(\Rset^{(n)} \times 0) \subset RC_n.
$
In fact, our approach to the Main Theorem is to prove the representation theoretic 
analogue of \eqref{Up1} for generic fibres and then to deduce the general case by continuity. In particular this avoids dealing directly with special fibres of $\Upsilon$,
such as those required for the Shapiro Conjecture.

Tracing back the relations between real loci and real representations,
we find interesting reformulations of the Shapiro conjecture in 
different setups. A curious one arises from the setting of the Wilson 
adelic Grassmannian:

\medskip
{\em{Fix a partition $\la= (\la_1 \geq \la_2 \geq \ldots \geq \la_l >0)$
and consider the corresponding Schur function $s_\la(p_1, p_2, \ldots, p_N)$,
$N = \la_1 + l -1$. If $c_1, c_2, \ldots c_N \in \Cset$ are such that 
$s_\la(x+ c_1, c_2, \ldots, c_N)$ has only real roots,
then $c_1, c_2, \ldots, c_N \in \Rset$.}}

\medskip
Finally we would like to point out that the real loci of other quiver varieties, 
and the reality of representations of other deformed preprojective
algebras of quivers and Cherednik algebras could be naturally 
related to other combinatorial problems of the Shapiro--Shapiro type.
\\ \hfill \\
%%%%%%%%%%%%%%%%%%%%%%%%%%%%%%%%%%%%%%%%%%%%%%%%%%%%%%%%%%%%%%%%%%%%
{\bf Acknowledgements.} We would like to thank Gwyn Bellamy, Yuri Berest,
Victor Ginzburg, Frank Sottile, and George Wilson for their 
helpful comments.
E.H. acknowledges the support by grant MI 1504/2005 of the National Fund 
"Scientific research" of the Bulgarian Ministry of Education and Science.
The research of M.Y. was supported by NSF
grant DMS-0406057 and an Alfred P. Sloan research fellowship.
%%%%%%%%%%%%%%%%%%%%%%%%%%%%%%%%%%%%%%%%%%%%%%%%%%%%%%%%%%%%%%%%%%%%%%%%%%
\sectionnew{Calogero-Moser spaces}
\label{CM}
First we recall the definition of Calogero--Moser spaces, for details 
we refer the reader to \cite{W2}. Define the locally closed 
subset $\ol{C}_n$ of $\gl_n(\Cset) \times \gl_n(\Cset)$, 
consisting of pairs of matrices $(X,Z)$ such that
\begin{equation}
    \rank ([X,Z] + I_n)=1,
\label{rk1}
\end{equation}
where $I_n$ is the identity matrix of size $n \times n$. The group $GL_n(\Cset)$ acts on $\ol{C}_n$ by simultaneous conjugation
\begin{equation}
          g.(X,Z) = (gXg^{-1}, gZg^{-1}), \,\, g \in GL_n(\Cset),
\end{equation}
and this action is free and proper, see \cite{W2}. There then exists
a geometric quotient
\begin{equation}
\pi_n \colon  \ol{C}_n \to C_n = \ol{C}_n/ GL_n(\Cset)
\end{equation}
which is a smooth, irreducible, complex affine variety, \cite{W2}.
It is called the $n$-th Calogero--Moser space. Define
the real locus of $\ol{C}_n$ by
\[
\ol{R}_n = \{(X,Z) \in \gl_n(\Rset) \times \gl_n(\Rset) \mid
\rank ([X,Z] + I_n)=1 \},
\]
and define the real locus of $C_n$ as the push-forward of the
real locus of $\ol{C}_n$ under $\pi_n$, namely
\[
RC_n := \pi_n \left( \ol{R}_n \right).
\]

The next proposition identifies $RC_n$ with 
a real Calogero--Moser space and explicitly describes it as 
a real algebraic subset of $C_n$.
Denote the natural inclusion
\[
\ol{i}_n \colon \ol{R}_n \hookrightarrow \ol{C}_n.
\]
Here we view $\ol{C}_n$ as a real variety and $\ol{i}_n$ 
as an embedding of real varieties. We also define
\begin{equation}
\O^\Rset(C_n) = \O(\ol{C}_n)^{GL_n(\Cset)} \cap 
\Rset[x_{jl}, z_{jl}]_{j,l=1}^n \subset \O(C_n),
\label{ORset}
\end{equation}
where $x_{jl}$, $z_{jl}$ are the matrix entries of $X,Z$ 
considered as regular functions on $\ol{C}_n$. 

\bpr{CMlocus} Keep the above notation. The set $RC_n$ coincides with the real 
algebraic subset of $C_n$ 
\[
\{ c \in C_n \mid f(c) \in \Rset \text{ for all }f \in \O^\Rset(C_n) \}.
\]

There exists a (smooth) geometric quotient for the 
action of $GL_n(\Rset)$ on $\ol{R}_n$ and thus the natural inclusion 
$\ol{i}_n \colon \ol{R}_n \hookrightarrow \ol{C}_n$ induces 
a morphism of real varieties 
$i_n \colon \ol{R}_n/GL_n(\Rset) \to C_n$. Furthermore,
$i_n$ defines an isomorphism of real varieties 
$i_n \colon \ol{R}_n/GL_n(\Rset) \cong RC_n$.
\epr

We will call the quotient $R_n := \ol{R}_n/GL_n(\Rset)$
the $n$-th real Calogero--Moser space.

\begin{proof} 
Consider the categorical quotient
$\nu_n \colon \ol{R}_n \to \ol{R}_n
/ \hspace{-0.08cm} / GL_n(\Rset)$. 
By the universal property of 
categorical quotients $\ol{i}_n$ descends to the map
(of real varieties)
\[
i_n \colon \ol{R}_n 
/ \hspace{-0.08cm} /
GL_n(\Rset) \to C_n.
\] 
The fiber of 
$\pi_n \ol{i}_n \colon \ol{R}_n \to C_n$ 
through $(X,Z) \in \ol{R}_n$ is 
\[
\{ (X_1, Z_1) \in \ol{R}_n \mid
\; \exists g \in GL_n(\Cset) \; 
\mbox{such that} \; X_1 = g X g^{-1}, 
Z_1 = g Z g^{-1} \}.
\]
\leref{conj} below implies that the fiber of
$i_n \nu_n = \pi_n \ol{i}_n \colon \ol{R}_n \to C_n$ 
through $(X,Z) \in \ol{R}_n$ is 
\[
\{ (X_1, Z_1) \in \ol{R}_n \mid
\; \exists g \in GL_n(\Rset) \; 
\mbox{such that} \; X_1 = g X g^{-1}, 
Z_1 = g Z g^{-1} \}.
\]
This implies that $i_n$ is an injection and that 
the fibers of $\nu_n$ are exactly the $GL_n(\Rset)$
orbits on $\ol{R}_n$. In particular, 
$\nu_n$ is a geometric quotient; hence 
we will denote its range by $\ol{R}_n/GL_n(\Rset).$
Analogously to \cite[Section 1]{W2} the latter
is smooth. 

It is clear that
\[
\Im i_n \subseteq RC_n \subseteq X_n :=
\{ c \in C_n \mid f(c) \in \Rset \text{ for all }f \in \O^\Rset(C_n) \}.
\]
Furthermore, $\O(C_n)$ is the complexification of $\O^\Rset(C_n)$
since as a $GL_n(\Rset)$ module 
$\O(\ol{C}_n) \cap \Rset[x_{jl}, z_{jl}]_{j,l=1}^n$ is a direct sum 
of finite dimensional modules and 
$\left(\O(\ol{C}_n) \cap \Rset[x_{jl}, z_{jl}]_{j,l=1}^n \right)_\Cset\cong
\O(\ol{C}_n)$. Thus
\[
i_n^* \colon \O(X_n) \to \O(R_n)
\]
is an isomorphism, where $\O(X_n)$ and $\O(R_n)$ denote the 
real coordinate rings of $X_n$ and $R_n$. This implies that 
$RC_n = X_n$ and that $i_n \colon R_n \to RC_n = X_n$ is an isomorphism 
of real algebraic varieties.
\end{proof}

\ble{conj} Let $X_1, \ldots, X_k, Y_1, \ldots, Y_k \in 
\gl_n(\Rset)$. If there exists $g \in GL_n(\Cset)$
such that
\begin{equation}
Y_1 = g X_1 g^{-1}, \ldots, Y_k = g X_k g^{-1}, 
\label{conj1}
\end{equation}
then there exists $g \in GL_n(\Rset)$ with the same property.
\ele
\begin{proof} We can assume that the element $g \in GL_n(\Cset)$
satisfying \eqref{conj1} is such that $\Re g$ is nondegenerate.
If this is not the case we can substitute $g$ with $a.g$ for an 
appropriate scalar $a \in \Cset^*$. Taking real parts 
in $Y_j g = g X_j$ we find
\[
Y_1 = (\Re g) X_1 (\Re g)^{-1}, \ldots, 
Y_k = (\Re g) Y_1 (\Re g)^{-1}, \quad \Re g \in GL_n(\Rset)
\] 
and thus $\Re g \in GL_n(\Rset)$ has the needed property.
\end{proof}

Recall the definition \eqref{map_Ups} of the finite map 
$\Upsilon \colon C_n \to \Cset^{(n)} \times \Cset^{(n)}$.
The following Theorem relates the real loci of the Calogero--Moser space
$C_n$ and $\Cset^{(n)} \times \Cset^{(n)}$ 
by the map $\Upsilon$. It was previously proved by 
Mukhin, Tarasov, and Varchenko \cite{MTV2}.

\bth{2} We have 
\[
\Upsilon^{-1}(\Rset^{(n)} \times \Rset^{(n)}) \subset RC_n.
\]
\ethe

We postpone the proof of \thref{2} to Sect. \ref{Cherednik}.

Denote by $\Cr$ the subset of $\Cset^n$ consisting 
of $(\la_1, \ldots, \la_n)$, $\la_j \neq \la_l$, for $j \neq l$, 
and by $\CrCS$ its image in $\Cset^{(n)}= \Cset^n/S_n$.
Wilson proved \cite{W2} that $C_n$ has a Zariski open subset 
isomorphic to $T^* \CrCS$. 
It is the image under $\pi_n$ of the subset of $\ol{C}_n$ 
consisting of pairs of matrices
\begin{equation}
X=\diag(\la_1, \ldots, \la_n), Z=                     
   \left(
   \begin{array}{cccc}
      \al_1 & (\la_1-\la_2)^{-1}  &\dots & (\la_1-\la_n)^{-1} \\
      (\la_2-\la_1)^{-1} & \al_2 &\dots  &(\la_2-\la_n)^{-1}  \\
      \dots                &  \dots    & \dots         & \dots  \\
      (\la_n-\la_1)^{-1} & (\la_n-\la_2)^{-1} & \dots & \al_n
   \end{array}
   \right)                     
\label{C-M} 
\end{equation}
where $(\la_1, \ldots \la_n) \in \Cr$, 
$(\al_1, \ldots, \al_n) \in \Cset^n$. 

The restriction of \thref{2} to this subset leads to the following 
Corollary for eigenvalues of Calogero--Moser matrices, previously 
proved by Mukhin, Tarasov, and Varchenko \cite{MTV2}.

\bco{1} Let $(\la_1, \la_2, \dots, \la_n)$ 
be an $n$-tuple of distinct real numbers. 
If an $n$-tuple $\alpha=(\alpha_1, \alpha_2,
\dots, \alpha_n) \in \Cset^n$ has the property that the 
Calogero--Moser matrix $Z$ in \eqref{C-M} has only real 
eigenvalues, then $\al \in \Rset^n$.
\eco

\begin{proof} Let $\la$ and $\mu$ 
be two real $n$-tuples as in \coref{1} and $X$ and $Z$ be the matrices \eqref{C-M}. 
Then according to \thref{2} there exists $g \in GL_n(\Cset)$, such
that $g.(X,Z) \in \gl_n(\Rset)^{\times 2}$. Since the
eigenvalues $\lambda_1, \dots, \lambda_n$ of the  matrix
$X$ are distinct, there exists a matrix $g_1\in GL_n(\Rset)$
for which  $g=g_1A,\,\, \textrm{where}\,\,A= diag(a_1,\dots,a_n)$
with some complex numbers $a_1, \dots, a_n$. The fact
that the matrix $gZg^{-1} = g_1AZA^{-1}g_1^{-1}$ is real shows that
$AZA^{-1}$ is also real. The off-diagonal entries of this
matrix are $(\lambda_i-\lambda_j)^{-1}a_ia^{-1}_j$, i.e.
$a_i a_1^{-1} \in \Rset$. One can write 
$(a_1,\dots, a_n) = a_1(b_1,\dots,b_n)$ 
for some real $n$-tuple $(b_1,\dots,b_n) \in \Rset^n$. 
Now put $g_2= g_1. \diag(b_1,\dots,b_n)$. Then 
$g_2 Z g_2^{-1} = g_1Zg_1^{-1}$ is  real and the matrix $g_2$ is real.
Thus $z$ has real entries and in particular $\al \in \Rset^n$.
\end{proof}
%%%%%%%%%%%%%%%%%%%%%%%%%%%%%%%%%%%%%%%%%%%%%%%%%%%%%%%%%%%%%%%%%%%%%%%%%%
\sectionnew{Cherednik algebras}
\label{Cherednik} Rational Cherednik algebras are two step
degenerations of double affine Hecke algebras \cite{Ch}. For
details on rational Cherednik algebras and more generally on
symplectic reflection algebras we refer the reader to the
Etingof--Ginzburg paper \cite{EG}. The rational Cherednik algebra
$H_{0,c}(S_n)$ is generated by the polynomial subalgebras
$\Cset[x_1, \ldots, x_n]$ and $\Cset[y_1, \ldots, y_n]$, and the
group algebra $\Cset S_n$ of the symmetric group, subject to the
following deformed crossed product relations
\begin{align}
&s_{ij} x_i = x_j s_{ij}, \quad s_{ij} y_i = y_j s_{ij},
\nn
\\
&[x_i, y_j]= c s_{ij} \; \; (i \neq j), \quad
[x_k, y_k] = - c \sum_{i \neq k} s_{ik}.
\nn
\end{align}
The algebras $H_{0,c}(S_n)$ are isomorphic for different values of
$c \neq 0$ and we will mostly restrict our attention to $c=1$.
Denote by $e= (1/n!) \sum_{\sigma \in S_n} \sigma$ the
symmetrizing idempotent of $\Cset S_n \subset H_{0,1}(S_n)$. The
spherical subalgebra of $H_{0,1}(S_n)$ is the subalgebra $U = e
H_{0,1}(S_n) e$, \cite{EG}.

First we recall several results of Etingof and Ginzburg on finite
dimensional irreducible $H_{0,1}(S_n)$ representations.

\bth{EG1} {\em{(}}Etingof-Ginzburg,
\cite[Theorems 1.23 and 1.24]{EG}{\em{)}}

(a) $U$ is a commutative algebra and is isomorphic to
the coordinate ring $\O(C_n)$ of the $n$-th Calogero--Moser
space $C_n$.

(b) The irreducible $H_{0,1}(S_n)$-representations are
parametrized by the points of $C_n$. Given $p \in C_n$, the
corresponding $H_{0,1}(S_n)$ representation is
\[
M_p := H_{0,1}(S_n) e \otimes_U \chi_p
\]
where $\chi_p \colon U \cong \O(C_n) \to \Cset$
is the character associated to $p$.

(c) Each representation $M_c$ has dimension $n!$ and, as an $S_n$
representation, is isomorphic to the regular representation of
$S_n$. \ethe

We will also need the following additional fact from \cite{EG}
regarding the structure of the representations $M_c$. First,
denote by $S_{n-1}$ the subgroup of $S_n$ permuting the last $n-1$
indices $\{2, \ldots, n\}$ and by $\ol{e} = (1/(n-1)!)
\sum_{\sigma \in S_{n-1}} \sigma$ the symmetrizing idempotent of
$S_{n-1}$. The subspace $\ol{e} M_c$ is stable under the action of
$x_1$ and $y_1$ because $x_1$ and $y_1$ commute with $\ol{e}$. In
any basis of $\ol{e} M$, $x_1$ and $y_1$ act by a pair of matrices
$(X_c, Z_c) \in \ol{C}_n$ such that $\pi_n(X_c, Z_c)=c$

Finally, for $c\in \Rset$ we denote the real subalgebra of
$H_{0,c}(S_n)$ generated by $x_1, \ldots, x_n,$ $y_1, \ldots, y_n$
and the elements of $S_n$ by $H_{0,c}^\Rset(S_n)$. It is clear
that $H_{0,c}(S_n)$ is the complexification of
$H_{0,c}^\Rset(S_n)$.

The following theorem is our main result.

\bth{Cher} Fix an irreducible $H_{0,1}(S_n)$ module $V$.
If the restriction of the operators $x_1$ and $y_1$ to $\ol{e} V$ 
have only real eigenvalues, then $V$ is the complexification 
of a (real) $H_{0,1}^\Rset(S_n)$ module.
\ethe

Before we prove \thref{Cher}, we note several lemmas. 
Let $\Rr= \Cr \cap \Rset^n$ and $\RrCS = \CrCS \cap \Rset^{(n)}$.

\ble{equiv}
(a) The statements of \thref{Cher} and \thref{2}
are equivalent.

(b) To prove \thref{2}, it is sufficient to show that 
\begin{equation}
\Y^{-1} (\RrCS \times \RrCS) \subset RC_n.
\label{reg1}
\end{equation}
To prove \thref{Cher}, it is sufficient to establish the validity 
of the statement for representations $V$ of $H_{0,1}(S_n)$ 
for which $x_1$ and $y_1$ act on $\ol{e} V$
by regular semisimple operators. 
\ele

\begin{proof} (a) Assume the validity of the statement
of \thref{Cher} and fix $(X,Z) \in \ol{C}_n$ such that both $X$
and $Z$ have only real eigenvalues. Set $p = \pi_n(X,Z)$.
\thref{Cher} implies that $M_p$ is the complexification of a real
$H_{0,1}^\Rset(S_n)$ representation $M^\Rset_p$. Therefore $\ol{e}
M_p^\Rset$ is a real vector space which is stable under $x_1$ and
$y_1$ and such that $(\ol{e} M_p^\Rset)_\Cset = \ol{e} M_p$. If
$X_1$ and $Z_1$ are the matrix representation of the restriction
of the operators $x_1$ and $y_1$ to $\ol{e} M_p^\Rset$ in any
basis of $\ol{e} M_p^\Rset$, then $(X_1, Z_1) \in \ol{R}_n \cap
GL_n(\Cset)(X,Z)$. Thus $\pi_n(X,Z) \in RC_n$.

In the opposite direction, let us assume the validity of the
statement of \thref{2}. Fix $p \in C_n$ such that the restriction
of $x_1$ and $y_1$ to $\ol{e} M_c$ have only real eigenvalues.
Then $p$ belongs to the real locus $RC_n$ of $C_n$. One checks by
a direct computation that under the Etingof--Ginzburg isomorphism
$U \cong \O(C_n)$, $U^\Rset = e H_{0,1}^\Rset(S_n) e$ corresponds
to $\O^\Rset(C_n)$, with the notation as in \eqref{ORset}. Because
of \prref{CMlocus}, $\chi_p \colon U \cong \O(C_n) \to \Cset$
restricts to a real character $\chi_p \colon U^\Rset \to \Rset$.
Then $M_p$ is the complexification of the $H_{0,1}^\Rset(S_n)$
representation $M_p^\Rset= H_{0,1}^\Rset(S_n) e \otimes_{U^\Rset}
\chi_p$, which establishes \thref{Cher}.

(b) Tracing back the equivalence in part (a), one sees that
the second statement is a consequence of the first one.
We proceed with the proof of the first statement.
Since the map $\Upsilon$ is open in the usual topology,
\[
\overline{\Upsilon^{-1}(S)} = \Upsilon^{-1}(\overline{S})
\]
for any subset $S$ of $\Cset^{(n)} \times \Cset^{(n)}$. Here
$\overline{(.)}$ refers to the closure in the usual topology. To
obtain the validity of \thref{2} from \eqref{reg1}, we apply this
for the set $S=\RrCS \times
\RrCS$ which is dense in
\[
{\Rset}^{(n)} \times {\Rset}^{(n)}.
\]
\end{proof}

We will need a very simple deformation theoretic argument for the
proof of \thref{Cher}.
%\ff{I had to add this because there was an oversight in my original 
%proof. All the evalues for $y_1$ are regular when we apply $\Theta$ to 
%them they do not necessarily end up in $\Rset[\Rr]^{S_n}$ but could 
%rather be smeared across $\Rset[\Rr\times \Rset^n]^{S_n}$. This lemma 
%allows us to use a little deformation argument to finesse this.}

\ble{deform}Suppose that $A$ is a commutative flat finite $\Rset
[t]$-algebra such that for $p\in \Rset$ each specialisation $A(p)
:= A/(t-p)A$ is a product of fields. Then $A(p) \cong A(q)$ for
all $p,q \in \Rset$. \ele
\begin{proof} We are going to prove that for any $p\in \Rset$ there is an open (analytic) interval $I$ containing $p$ such that $A(v) \cong A(p)$ for any $v\in I$.

Let $f_1(p), \ldots , f_m(p)$ be a complete set of primitive
idempotents for $A(p)$ that extends to an $\Rset$-basis of $A(p)$,
$f_1(p), \ldots , f_m(p), f_{m+1}(p), \ldots , f_n(p)$, such that
$f_i(p)f_{m+j}(p) = \delta_{ij}f_{m+j}(p)$ for $1\leq i\leq m$ and
$1\leq j \leq n-m$. Since $A$ is a free $\Rset[t]$-module, this
lifts to a basis of $A$ denoted $f_1, \ldots , f_n$. For any $u\in
\Rset$ we will denote the induced multiplication in $A(u)$ by
$\ast_u$ so that $f_i(u) \ast_u f_j(u) = \sum_{k=1}^n
\alpha_{i,j}^k(u) f_k(u)$ where $\alpha_{i,j}^k \in \Rset[t]$. We
have
\begin{equation} \label{strcons} \alpha_{i,j}^k (p) =
\delta_{i,j}\delta_{i,k} \text{ for $1\leq i,j \leq
m$}\end{equation} and \begin{equation} \label{strcons2}
\alpha_{i,m+j}^k (p) = \delta_{i,j}\delta_{m+j, k} \text{ for
$1\leq i\leq m$, $1\leq j \leq n-m$}.
\end{equation}

Consider the function $G: \Rset\times \Rset^n \longrightarrow
\Rset^n$ that sends $(u; \lambda_1, \ldots , \lambda_n)$ to the
coefficients of $(\sum \lambda_i f_i(u)) \ast_u (\sum \lambda_i
f_i(u)) - (\sum \lambda_i f_i(u))$  in the basis $(f_1(u), \ldots
, f_n(u))$ of $A(t)$. In other words $G = (G_1, \ldots , G_n)$
where
$$G_k (t; \lambda_1, \ldots , \lambda_n) = \sum_{i,j}
\lambda_i\lambda_j \alpha_{i,j}^k- \lambda_k.$$ For $1\leq i \leq
m$ set $p_i =(p; 0, \ldots ,0,1,0, \ldots 0)$, the $1$ occurring in
the $i$th place. Then, by \eqref{strcons} and \eqref{strcons2}, $
(
\partial{G_k}/\partial{\lambda_l})_{k,l}(p_i)$ is a diagonal
matrix with entries from $\{1,-1\}$. Thus the determinant at $p_i$
is non-zero, so there exists an open interval $U_i$ containing $p$
and continuous mapping $\theta_i  : U_i \longrightarrow \Rset^n$
such that $\tilde{F}_i(u) = \sum_{j=1}^n (\theta_i(u))_j  f_j(u)$
is an idempotent for each $u\in U_i$ and $\tilde{F}_i(p) =
f_i(p)$. Repeating this argument for each $1\leq i \leq m$
produces a continuous family of idempotents $\tilde{F}_1(u),
\ldots , \tilde{F}_m(u)$ for every $u\in U=\bigcap_i U_i$ and such
that $\tilde{F}_i(p) = f_i(p)$. This family can be adjusted
inductively to produce orthogonal idempotents $F_1(u), \ldots ,
F_m(u)$. Indeed we set $F_1(u) = \tilde{F}_1(u)$. Then if we have
found $\tilde{F}_1(u), \ldots , \tilde{F}_{s-1}(u)$ for some $s<
m$ we set $F_{s}(u) = (1-F_1(u)- \cdots -
F_{s-1}(u))\tilde{F}_{s}(u)$. We finish by setting $F_m(u) = 1-
F_1(u) - \cdots - F_{m-1}(u)$.

Now for any $u\in U$ we have an algebra decomposition $$A(u) =
\bigoplus_{j=1}^m F_j(u)A(u)F_j(u)= \bigoplus_{j=1}^m A(u)_j.$$
Let $1\leq j \leq n-m$ so that $\dim A(p)_j =2$  with basis
$f_j(p), f_{m+j}(p)$. Then on some open set $V_j$ of $U$ we see
that $F_{m+j}(u):= F_j(u) f_{m+j}(u)$ is non-vanishing, so that
$\dim A(v)_j \geq 2$ for all $v\in V_j$. Hence we produce an open
interval $V = U\cap (\bigcap_{j=1}^{n-m} V_j)$ which contains $p$
and on which $\dim A(v)_i \geq \dim A(p)_i$ for $v\in V$ and all
$1\leq i \leq m$. Since $\dim A(v) = \sum \dim A(v)_i$ is constant
we find $\dim A(v)_i = \dim A(p)_i$ for all $1\leq i \leq m$ and
$v\in V$.

We must show that $A(v)_i \cong A(p)_i$. This is obvious if
$\dim_{\Rset} A(v)_i = 1$, so we assume that $\dim A(v)_i = 2$. By
the previous paragraph we have a basis $F_i(v), F_{m+i}(v)$ for
$A(v)_i$ where $F_i(v)$ is the identity element of $A(v)_i$. Since
$A(p)_i \cong \Cset$ we may assume without loss of generality that
$F_{m+i}(p)\ast_p F_{m+i}(p) = -F_i(p)$, so if we write
$F_{m+i}(v)\ast_v F_{m+i}(v) = \alpha(v) F_i(v) + \beta(v)
F_{m+i}(v)$, then $\alpha(p) = -1$ and $\beta (p) = 0$. Now
consider the mapping $H: V \times \Rset^2 \longrightarrow \Rset^2$
given by sending $(v ; \lambda_1 , \lambda_2 )$ to the
coefficients in the basis $F_i(v), F_{m+i}(v)$ of the expression
$$(\lambda_1 F_i(v) + \lambda_2 F_{m+i}(v) )\ast_v (\lambda_1
F_i(v) + \lambda_2 F_{m+i}(v)) + F_i(v).$$ In other words $H(v;
\lambda_1, \lambda_2 ) = (H_1, H_2) = (\lambda_1^2 + 1 +
\lambda_2^2\alpha(v), 2\lambda_1\lambda_2 + \lambda_2^2\beta(v)).$
At $(p; 0,1)$ $\det (\partial H_i /\partial \lambda_j) = 4$, so we
can find an open interval $I\subseteq V$ including $p$ and for all
$v\in I$ a basis $F_i(v), X_i(v)$ of $A(v)_i$ such that
$X_i(v)\ast_v X_i(v) = -F_1(v)$. Hence $A(v)_i \cong \Cset \cong
A(p)_i$, as required.
\end{proof}

The lemma that follows is well-known, but we have been unable to
find a proof in the literature. To ease notation let $\alpha_{ij}
= x_i - x_j$ and $\alpha_{ij}^{\vee} = y_i-y_j$ and let $\delta =
\prod_{i<j} \alpha_{ij}\in \Rset[x_1, \ldots , x_n]$ be the
discriminant. 

\ble{dunkemb} For any $c\in \Cset$ we have $H_{0,c}[\delta^{-1}]
\cong \Cset[\Cr\times \Cset^n] \ast S_n$ and if $c\in \Rset$ then
$H_{0,c}^{\Rset}[\delta^{-1}] \cong \Rset[\Rr\times \Rset^n] \ast
S_n$. \ele
\begin{proof}
Assume that $c\in \Rset$. We define the isomorphism
$\Theta_c^{\Rset}: H_{0,c}^{\Rset}[\delta^{-1}] \longrightarrow
\Rset[\Rr\times \Rset^n] \ast S_n$ by $$x \mapsto x, \quad y
\mapsto y + c\Delta(y),  \quad w\mapsto w,$$ where $\Delta(y) =
\sum_{i<j} \frac{\langle y, \alpha_{ij}\rangle}{\alpha_{ij}}((i\,
j)-1)$. Once we show that this mapping is well-defined it is clear
that it is an isomorphism since we can remove the $\Delta(y)$-term
in $\Theta_c^{\Rset}(y)$ by subtracting elements from
$\Rset[\Rr]\ast S_n$.

To prove well-definedness we recall the Dunkl isomorphism
constructed in \cite[Proposition 4.5]{EG}. For $t\neq 0$ this
produces an isomorphism $\Theta_c^t: H_{t,c}[\delta^{-1}]
\longrightarrow \mathcal{D}(\Cr)\ast S_n$ by $$x \mapsto x, \quad
y \mapsto t\partial_y + c\sum_{i<j}  \frac{\langle y,
\alpha_{ij}\rangle}{\alpha_{ij}}((i\, j)-1), \quad w\mapsto w.$$
Thus $\Theta_c^{\Rset}$ is a real form of the semi-classical limit
$\mathop{\lim}_{t\rightarrow 0} \Theta_c^t$. From the Dunkl
isomorphism $\Theta_1^t$ and the defining relation for $H_{t,1}$, see \cite[Formula (1.15)]{EG}, we find
$$[\Delta(y), x] = -\frac{1}{2}\sum_{i<j } \langle y,
\alpha_{ij}\rangle \langle \alpha_{ij}^{\vee} , x\rangle(i\, j)
\text{ and } [\Delta(y), \Delta(y') ] = 0.$$ From this we deduce
that
$$[\Theta_c^{\Rset}(y), \Theta_c^{\Rset}(x) ] = [y,x] + c[\Delta(y),
x] = -c\frac{1}{2}\sum_{i<j } \langle y,\alpha_{ij}\rangle \langle
\alpha_{ij}^{\vee} , x\rangle(i\, j) = \Theta_c^{\Rset}([y,x]),$$
and \begin{eqnarray*} [\Theta_c^{\Rset}(y), \Theta_c^{\Rset}(y')] &=&[ y,y']+
c[y,\Delta(y')] + c[ \Delta(y), y'] + c^2[\Delta(y), \Delta(y')] \\ &=&
c([y,\Delta(y')] + [ \Delta(y), y']) .\end{eqnarray*} So it remains to check
that $[y,\Delta(y')] + [ \Delta(y), y'] =0$. Well,
\begin{eqnarray*}
[y,\Delta(y')] + [ \Delta(y), y'] &=& \sum_{i< j} \left( \frac{\langle y', \alpha_{ij} \rangle}{\alpha_{ij}} [y, (i\, j)-1] - \frac{\langle y,\alpha_{ij}\rangle}{\alpha_{ij}} [y', (i\,j)-1] \right) \\ &= & \sum_{i< j} \frac{1}{\alpha_{ij}}\left(\langle y', \alpha_{ij}\rangle (y - {}^{(i\,j)}y) -\langle y, \alpha_{ij} \rangle (y' - {}^{(i\, j)}y') \right)(i\, j) \\ & = &\sum_{i< j} \frac{1}{\alpha_{ij}}\left( \langle y', \alpha_{ij} \rangle  \langle y, \alpha_{ij} \rangle \alpha_{ij}^{\vee} - \langle y,\alpha_{ij} \rangle \langle y', \alpha_{ij} \rangle\alpha_{ij}^{\vee}  \right)(i \, j) \\ &=& 0.
\end{eqnarray*}
The isomorphism for $H_{0,c}[\delta^{-1}]$ follows by an identical
argument.
\end{proof}

%\begin{theorem}
% Let $V$ be an irreducible representation of $H$ such that the restriction of the operators $x_1$ and $y_1$ to $\overline{e}V$ have only real eigenvalues. Then $V = \Cset\otimes_{\Rset} V^{\Rset}$ where $V^{\Rset}$ is an irreducible $H^{\Rset}$-module.
%\end{theorem}
\noindent {\it Proof of \thref{Cher}.} By  \leref{equiv} it is
enough to prove this when both $x_1$ and $y_1$ have distinct
eigenvalues on $\overline{e}V$. Now by \cite[Theorem 11.16]{EG} the
eigenvalues of $x_1$ on $\overline{e}V$ coincide with the action
of $\Cset[\Cset^n]^{S_n} \subset Z(H)$ on $V$ and so we see that
$\delta^2 \in \Cset[\Cset^n]^{S_n}$ acts by a non-zero scalar on
$V$. In particular $V$ is naturally an irreducible
$H[\delta^{-1}]$-representation.

Let $\mathfrak{m}$ be the maximal ideal of
$\Rset[\Rr]^{S_n}\otimes \Rset[\Rr]^{S_n}$ corresponding to the
eigenvalues of $x_1$ and $y_1$. For $c\in \Rset$ set $$H(c) =
H^{\Rset}_{0,c}/\mathfrak{m}H^{\Rset}_{0,c} \cong
H^{\Rset}_{0,c}[\delta^{-1}]/\mathfrak{m}H^{\Rset}_{0,c}[\delta^{-1}],$$
a flat family of algebras over $\Rset$. By definition $V$ is an
irreducible $\Cset\otimes_{\Rset} H(1)$-representation, so we must
prove that all such representations are extensions of
$H(1)$-representations.
%So by the above lemma we are reduced to studying the following situation: $V$ is an irreducible $\Cset[\Cr\times \Cr]\ast S_n$-representation such that the scalar function $\chi_V$ by which the central subalgebra $\Cset[\Cr]^{S_n}\otimes \Cset[\Cr]^{S_n}$ acts is real-valued. Let $\mathfrak{m} \triangleleft \Rset[\Rr]^{S_n}\otimes \Rset[\Rr]^{S_n}$ be the corresponding maximal ideal associated to $\chi_V$, so that $V$ is an irreducible $\overline{H}$-module where $$\overline{H} = \frac{\Cset[\Cr\times \Cr]}{\mathfrak{m}\Cset[\Cr\times \Cr]}\ast S_n \cong \Cset\otimes_{\Rset} \overline{H}^{\Rset}$$ where $\overline{H}^{\Rset} = \left( \frac{\Rset[\Rr\times \Rr]}{\mathfrak{m}\Rset[\Rr\times \Rr]}\ast S_n\right).$
To do this we will simply prove that $H(c) \cong
\Mat_{n!}(\Rset)^{\oplus n!}.$

We translate the problem to $\Rset[\Rr\times \Rset^n]\ast S_n$ by
applying $\Theta_c^{\Rset}$ of \leref{dunkemb} to the family
$H(c)$. This produces the algebras $$\tilde{H}(c) =
\frac{\Rset[\Rr\times \Rset^n]\ast
S_n}{\Theta_c^{\Rset}(\mathfrak{m})\Rset[\Rr\times \Rset^n]\ast S_n}$$
However, $\Theta^{\Rset}_c(\mathfrak{m}) \subset Z(\Rset[\Rr\times
\Rset^n]\ast S_n) = \Rset[\Rr\times \Rset^n]^{S_n}$ so that we actually
have $$\tilde{H}(c) = \left( \frac{\Rset[\Rr\times
\Rset^n]}{\Theta_c^{\Rset}(\mathfrak{m})\Rset[\Rr\times
\Rset^n]}\right)\ast S_n.$$ Hence we find a flat family of commutative
$\Rset$-algebras of dimension $(n!)^2$ $$A(c) =
\frac{\Rset[\Rr\times
\Rset^n]}{\Theta_c^{\Rset}(\mathfrak{m})\Rset[\Rr\times \Rset^n]}.$$

We consider first $c=0$. The homomorphism $\Theta_0$ restricts to
the inclusion $\Rset[\Rset^n]^{S_n}\otimes \Rset[\Rset^n]^{S_n}
\longrightarrow \Rset[\Rr\times \Rset^n]$ and so $$A(0)
=\frac{\Rset[\Rr\times \Rset^n]}{\mathfrak{m}\Rset[\Rr\times \Rset^n]} \cong \frac{\Rset[\Rr\times \Rr]}{\mathfrak{m}\Rset[\Rr\times \Rr]},$$ where the last isomorphism holds since $\mathfrak{m}$ is a maximal ideal of $\Rset[\Rr]^{S_n}\otimes \Rset[\Rr]^{S_n}$.
Since $S_n\times S_n$ acts freely on $\Rr\times \Rr$ it follows
that there are exactly $(n!)^2$ points of $\Rr\times \Rr$ lying
above $\mathfrak{m}$ and, in particular, $n!$ distinct free
(diagonal) $S_n$-orbits of points in $\Rr\times \Rr$. Thus we have
a direct product expansion $$A(0) \cong \prod_{i=1}^{n!}
\prod_{\sigma\in S_n} \Rset e_{i, \sigma}$$ where the
$e_{i,\sigma}$ are pairwise orthogonal idempotents and $\tau
e_{i,\sigma} = e_{i, \tau\sigma}$ for any $\tau\in S_n$.

Note that each $A(c)$ is separable. To see this we may assume that $c\neq 0$ since we dealt with the case $c=0$ above. Consider the
largest nilpotent ideal $J$ of $A(c)$. It must be $S_n$-stable and
so extend to a nilpotent ideal $(J\otimes \Cset)\ast S_n$ of
$\tilde{H}(c)\otimes \Cset \cong H(c)\otimes \Cset$. But this last
algebra is semisimple by \cite[Theorem 1.7(i)]{EG} and so $J=0$.

Now \leref{deform} shows that all $A(c)$ are isomorphic as
$\Rset$-algebras. Thus since the maximal spectrum of $A(0)$
corresponds to $n!$ distinct free (diagonal) $S_n$-orbits in
$\Rr\times \Rr$, the same is true for all $A(c)$. Therefore
${H}(c)\cong \prod_{i=1}^{n!} \left( \prod_{\sigma \in S_n}
e_{i,\sigma} \right)\ast S_n.$ But $\left(\prod_{\sigma \in S_n}
e_{i,\sigma} \right)\ast S_n \cong \Mat_{n!}(\Rset)$, the
isomorphism being given by sending $e_{i,\sigma}\otimes \tau$ to
the elementary matrix $E_{\sigma, \tau^{-1}\sigma}$, where we
label the rows and columns of $\Mat_{n!}(\Rset)$ by the elements of
$S_n$. Setting $c=1$ proves the theorem. $\hfill \Box$

Note that \thref{2} follows from \thref{Cher} and \leref{equiv} (a).
%%%%%%%%%%%%%%%%%%%%%%%%%%%%%%%%%%%%%%%%%%%%%%%%%%%%%%%%%%%%%%%
\sectionnew{Deformed preprojective algebras of quivers}
\label{quiver}
Let $Q=(Q_0, Q_1)$ be a finite quiver with vertex set $Q_0$ and 
arrow set $Q_1$. Denote by $\ol{Q}$ its double, obtained
by adding a reverse arrow $a^*$ for each arrow $a$ of $Q_1$.
The deformed preprojective algebra of $Q$ of weight $\nu=(\nu_i)_{i \in Q_0} \in \Cset^{Q_0}$, 
was defined by Crawley-Boevey and Holland \cite{CBH} 
as the quotient of the path algebra $\Cset\ol{Q}$
\[
\Pi_{\nu} (Q) = \Cset \ol{Q} / \langle \sum_{a \in Q_1} [a,a^*] - 
\sum_{i \in Q_0} \nu_i e_i \rangle,
\] 
where the $e_i$ denote the standard idempotents of 
$\Cset\ol{Q}$. For a real 
weight $\nu=(\nu_i)_{i \in I} \in \Rset^I$, let
\[
\Pi_{\nu}^\Rset (Q) = \Rset \ol{Q} / \langle \sum_{a \in Q_1} [a,a^*] - 
\sum_{i \in Q_0} \nu_i e_i \rangle.
\] 
Thus, for $\nu \in \Rset^I$, 
$\Pi_{\nu} (Q)$ is the complexification of 
$\Pi_{\nu}^\Rset (Q)$.

We restrict our attention to the quiver $Q$ with 2 vertices
$0$ and $\infty$, and two arrows $v \colon 0 \to \infty$ 
and $X \colon 0 \to 0$. We set $\nu_0=-1, \nu_\infty= n$ and denote $w= v^*$, $Z= X^*$. The algebra 
$\Pi_{\nu}(Q)$ is then generated 
by $X,Z,v,w$ and the idempotents $e_0, e_\infty$, and these satisfy the path algebra relations and 
\[
[X,Z] - wv = - e_0, \quad vw = n e_\infty.
\]
A left $\Pi_{\nu}(Q)$-module is thus a complex vector 
space $V = V_0 \oplus V_\infty$ with the data of 
\[
X, Z \in \End (V_0), \quad v \in \Hom (V_0, V_\infty), 
\quad w \in \Hom(V_\infty, V_0)
\]
such that 
\[
[X,Z] + \Id_{V_0} = w v, \quad 
vw = n \Id_{V_\infty}.
\] 
We restrict our attention to representations 
of $\Pi_{\nu}(Q)$ of dimension vector $(n,1)$, i.e. such that $\dim V_0 = n$, $\dim V_\infty =1$.
All such representations are irreducible by \cite{CB, W2}.

\thref{2} has the following corollary.

\bco{4} Consider a representation $V_0 \oplus V_\infty$ 
of the deformed preprojective 
algebra $\Pi_{\nu}(Q)$ of the above quiver for the  
weight $\nu= (-1, n)$ with dimension vector $(n,1)$.
If the operators $X, Z \in \End(V_0)$ have only real eigenvalues,
then this representation is the complexification of a 
representation of the real algebra $\Pi_{\nu}^\Rset(Q)$.
\eco

\begin{proof}
By \thref{2} we can find a basis of $V_0$ for which the entries of $X$ and $Z$ are real. Then any non-zero element in the image of this basis under the mapping $v$ provides a basis of $V_\infty$ for which the entries of $v$ and $w$ are real.
\end{proof}
%%%%%%%%%%%%%%%%%%%%%%%%%%%%%%%%%%%%%%%%%%%%%%%%%%%%%%%%%%%%
\sectionnew{The real locus of the Wilson's adelic Grassmannian}
\label{Wilson}
In this section we define the real locus of Wilson's 
adelic Grassmannian. Following \thref{2}
we formulate an elementary criterion for a point of Wilson's
Grassmannian to belong to its real locus. 

We start with a few general facts on the real locus of
Wilson's adelic Grassmannian. For details, we
refer the reader to Wilson's papers \cite{W1, W2}, 
van Moerbeke's review \cite{vM}, and the paper \cite{BHY2}. 
 
Sato's Grassmannian is an infinite dimensional Grassmannian 
of subspaces of $\Cset[z][[z^{-1}]]$ of a particular type, see 
\cite{vM}. It is the phase space of the KP hierarchy (an 
infinite dimensional integrable system). Wilson's adelic Grassmannian
is the subset of Sato's Grassmannian
which parametrizes all rank 1 bispectral wave 
functions, \cite{W1}. To a point $W \in \Gr^{\ad}$, in other words a subspace of 
$\Cset[z][[z^{-1}]]$ of a particular type, one associates its tau 
function
\[
\tau_W(t_1, t_2, \ldots) \in \Cset[[t_1, t_2, \ldots]]
\]
(defined up to a nonzero factor)
which is the image of $W$ under the Pl\"ucker 
embedding of $\Gr^{\ad}$ into the projectivization of 
$\Cset[[t_1, t_2, \ldots]]$. To $W \in \Gr^{\ad}$ one 
also associates its wave function
\begin{equation}
\Psi_W(x,z) = e^{xz}\left( 1 + a_1(x) z^{-1} + a_2(x) z^{-2} + \cdots 
\right), \; \; a_1(x) , a_2(x), \ldots \in \Cset(x).  
\label{Psi}
\end{equation}
The tau and wave functions of $W$ are related by Sato's 
formula
\begin{equation}
\Psi_W(x,z)= e^{\sum_{k=1}^\infty t_kz^k}
\frac{\tau\left(t-[z^{-1}]\right)}{\tau(t)} 
\Big|_{t_1=x, t_2 = t_3 = \ldots = 0}.
\label{Sato}
\end{equation}
Here and below we abbreviate $t = (t_1, t_2, \ldots)$ 
and $[z^{-1}]= (z^{-1}, z^{-2}/2, z^{-3}/3, \ldots)$. 

We will need two different parametrizations of $\Gr^{\ad}$, both due to 
Wilson, \cite{W1,W2}. According to \cite[Proposition 2.9]{W2} 
all tau functions $\tau_W(t)$ are polynomials 
in $t_1$ with leading coefficient 1 -- the other coefficients 
depend on $t_2, t_3, \ldots$. Denote by 
$\Gr^{\ad}_n$ the set of tau functions in $\Gr^{\ad}$ which 
are polynomials of degree $n$ in $t_1$.
The coefficients 
of $\tau_W$ give $\Gr^{\ad}_n$ the structure 
of a finite-dimensional complex affine variety. 
Set $C= \sqcup_{n \in \Zset_{\geq 0}} C_n$ and define Wilson's 
map $\beta \colon C \to \Gr^{\ad}$ by
\begin{equation}
\beta(\pi_n(X,Z))=W \;\; 
\mbox{where} \; \;  \tau_W:= 
\det (X+ \sum_{j=1}^{\infty}jt_j(-Z)^{i-1}), 
\; \mbox{for} \; (X,Z) \in \ol{C}_n.
\label{C-M-tau}
\end{equation}
Clearly $\beta(C_n) \subset \Gr^{\ad}_n$. The corresponding wave
function is given by
\begin{equation}
\Psi_W = e^{xz} \det ( I_n - (xI_n +X)^{-1} (zI_n + Z)^{-1} ),
\label{wave_beta}
\end{equation}
where $I_n$ is the identity matrix of size $n \times n$.

We denote by $G$ the set consisting of a finite
collection of distinct complex numbers $\mu_1, \ldots, \mu_k$
$(k \in \Zset_{\geq 0})$ and a collection of  
subspaces $V_1, \ldots, V_k$ of $\Cset[x]$, associated to each
of them. In other words $G$ consists of tuples
\[
(\mu_1, \ldots, \mu_k, V_1, \ldots, V_k)
\]
where two tuples of this kind are identified
if one of them is obtained from the other by 
a simultaneous permutation of the $\mu$'s and the $V$'s.
Define Wilson's map \cite{W1}
$\gamma \colon G \to \Gr^{\ad}$ by
\begin{equation}
\gamma(\mu_1, \ldots, \mu_k, V_1, \ldots, V_k)= W \; \; 
\mbox{where} \; \; \Psi_W(x,z) =
\frac{1}{p(z)} P_W(x, \dx) e^{xz}.
\label{gamma}
\end{equation}
Here
$p(z) = \prod_j (z - \mu_j)^{\dim V_j}$ and 
$P_W(x, \dx)$ is the monic differential operator
of degree $\dim V_1 + \cdots + \dim V_k$ and kernel
\[
e^{\mu_1 x} V_1 \oplus \cdots \oplus e^{\mu_k x} V_k.
\]
(This is the version of Wilson's map from \cite{BHY2}.)
Finally, denote by $G'$ the subset of $G$ consisting
of $(\mu_1, \ldots, \mu_k, V_1, \ldots, V_k)$
where all subspaces $V_j$ of $\Cset[x]$ contain no
constants except $0$. 

\bth{Wilson} {\em{(}}Wilson, \cite{W1,W2}{\em{)}}

(a) The maps $\beta \colon C_n \to \Gr^{\ad}_n$ are 
bijections.

(b) The map $\gamma \colon G' \to \Gr^{\ad}$ is a bijection.
Moreover, for a pair $(y \in G, y' \in G')$, 
$\gamma(y') = \gamma(y)$ if an only if $y$ is obtained from $y'$
by repeated applications of one of the following two rules:
\[
(\mu_1, \ldots, \mu_k, V_1, \ldots, V_k) \mapsto
(\mu_1, \ldots, \mu_k, \mu_{k+1}, V_1, \ldots, V_k, \Cset)
\]
(for $\mu_{k+1} \neq \mu_j$, $j = 1, \ldots, k$) and
\[
(\mu_1, \ldots, \mu_k, V_1, \ldots, V_k) \mapsto
(\mu_1, \ldots, \mu_k, V_1, \ldots, V_{k-1}, \wt{V}_k)
\]
where $\wt{V}_k = \{ p(x) \mid p'(x) \in V_k \}$.
\ethe

In the next theorem we describe the set of real points of 
Wilson's adelic Grassmannian.

\bth{App} For a point $W$ in Wilson's adelic Grassmannian $\Gr^{\ad}$
the following conditions are equivalent.

1) The plane $W$ has a real basis.

2) Up to a nonzero scalar the tau function $\tau_W$ has 
real coefficients.

3) The wave function $\Psi_W(x,z)$ has real coefficients, i.e.
$a_1(x), a_2(x), \ldots \in \Rset(x)$ in \eqref{Psi}.

4) $W= \beta(c)$ for some $n \in \Zset_{\geq 0}$, $c \in RC_n$,
i.e. $\tau_W$ is given by \eqref{C-M-tau} for a pair of real matrices 
$(X, Z) \in \ol{R}_n$.

5) $W= \gamma(\mu_1, \ldots, \mu_k, V_1, \ldots, V_k)$ for some 
$(\mu_1, \ldots, \mu_k, V_1, \ldots, V_k) \in G$ 
with the following properties: 

(a)  each $\mu_j$ appears together with its complex
conjugate and

(b) if $\ol{\mu}_j = \mu_k$, then $\ol{V}_j= V_k$. 
\ethe
Here for $p(x) = \sum_{j=1}^m a_j x^j$, we set
$\ol{p(x)} = \sum_{j=1}^m \ol{a}_j x^j$. Condition 
4 is the same as saying that {\em{one can find bases 
of $V_1, \ldots, V_k$ such that the basis of $V_j$ 
is real if $\mu_j$ is real and the bases of $V_j$ and $V_k$
are complex conjugate if $\ol{\mu}_j = \mu_k$}}. It can be also 
restated to: {\em{the space of functions 
$e^{\mu_1 x} V_1 \oplus \cdots \oplus e^{\mu_k x} V_k$ on 
$\Rset$ has a real basis}}.

We define the {\em real locus of Wilson's adelic Grassmannian}
as the set of all $W \in \Gr^{\ad}$ which satisfy any of the five 
equivalent conditions in \thref{App}.

\begin{proof} First we show 
$1 \Rightarrow 2 \Rightarrow 3 \Rightarrow 1$.
The tau function $\tau_W(t)$ is constructed from $W$ by the so called 
boson-fermion correspondence \cite{K}, from which $1 \Rightarrow 2$ is 
straightforward. Sato's formula \eqref{Sato} shows 
that $2 \Rightarrow 3.$
Finally $3 \Rightarrow 1$ because 
\[
W = \Span \{ \dx^j \Psi_W(x, z)|_{x=0} \}_{j=0}^\infty,
\]
cf. \cite{vM}. (If $\Psi_W(x,z)$ is singular at 0, one evaluates
its derivatives at any real point $r$. To get $W$, one flows back 
the plane obtained in this fashion, with respect to the first KP flow,
see e.g. \cite{BHY2}.)

The implications $4 \Rightarrow 2$ and $5 \Rightarrow 3$ are obvious
from \eqref{C-M-tau} and \eqref{gamma}.

Next we show $2 \Rightarrow 5$. Assume that up to a nonzero scalar
$\tau_W(t)$ has real coefficients and that 
$W = \gamma(\mu_1, \ldots, \mu_k, V_1, \ldots, V_k)$, for some
$(\mu_1, \ldots, \mu_k, V_1, \ldots, V_k) \in G'$. Taking complex
conjugates in \cite[Theorem 1]{BHY1} (see also \cite[eq. (5.7)]{W1}), 
we find $W = \gamma(\ol{\mu}_1, \ldots, \ol{\mu}_k, 
\ol{V}_1, \ldots, \ol{V}_k)$. Now 5 follows from  
the bijectivity of $\gamma \colon G' \to \Gr^{\ad}$, 
and the fact that $(\ol{\mu}_1, \ldots, \ol{\mu}_k, 
\ol{V}_1, \ldots, \ol{V}_k) \in G'$. We recall that two 
tuples of the type 
$(\mu_1, \ldots, \mu_k, V_1, \ldots, V_k)$ are identified
as elements of $G'$, if one of them is obtained from the 
other by simultaneous permutations of the $\mu$'s and the 
$V$'s.

Finally we prove $2 \Rightarrow 4$. We need to show that 
for $c \in C_n$, if $\tau_{\beta(c)}(t)$
has real coefficients up to a non-zero factor, then 
$c \in RC_n$. First we assume that $c = \pi_n(X,Z)$ for some
$(X,Z) \in \ol{C}_n$ where $X$ has distinct eigenvalues.
Since $X$ has distinct eigenvalues, we can assume that 
$X = \diag(\la_1, \ldots, \la_n)$ and that $Z$ has the form 
\eqref{C-M} for some $\la_1, \ldots, \la_n, \al_1, \ldots, \al_n 
\in \Cset$.
The reality of $\tau_{\beta(c)}(t)$ implies that
$\beta(\pi_n(X,Z)) = \beta(\pi_n(\ol{X}, \ol{Z}))$.
Because $\beta \colon C_n \to \Gr^{ad}_n$ is a bijection by
\thref{Wilson} (a), there exists $g \in GL_n(\Cset)$
such that $\ol{X} = g X g^{-1}$, $\ol{Z}= g Z g^{-1}$.
In particular, possibly after reindexing, 
we have for some $l \leq n/2$:
\[
\la_1, \ldots, \la_{2 l} \notin \Rset, 
\la_2 = \ol{\la}_1, \ldots, \la_{2l} =\ol{\la}_{2l-1}, 
\quad \la_{2l+1}, \ldots, \la_n \in \Rset
\]
and $\ol{Z} = g Z g^{-1}$ for
$
g = \diag(B, \ldots, B, I_{n-2l}), \; \; 
\mbox{where} \; \; 
B= 
\left(
      \begin{array}{cc}
         0& 1 \\
         1 & 0
      \end{array}
   \right).
$
This implies that 
\[
\al_2 = \ol{\al}_1, \ldots, \al_{2l} = \ol{\al}_{2l-1}, \quad
\al_{2l+1}, \ldots, \al_n \in \Rset.
\]
Set 
\[
B'= \left(
      \begin{array}{cc}
         1& 1 \\
         i &-i
      \end{array}
   \right), \; \; 
S = \diag( B', \ldots, B', I_{n-2l})
\]
and $X' = S X S^{-1}, Z' = S Z S^{-1}$.
Then the matrices $X'$ and $Z'$ have real entries and $c = \pi_n(X',Z')$. Thus $c \in RC_n$.

For the general case of $2 \Rightarrow 4$, suppose that, up to a nonzero factor, $\tau_{\beta(c)}(t)$ has real coefficients for some $c \in C_n$. 
Because of \cite[Proposition 8.6]{SW} for almost all 
$r \in \Rset$, $\tau_{\beta(c)}(t_1+r, t_2, \ldots) = 
\tau_{\beta(c')}(t_1, t_2, \ldots)$ for some $c' \in C_n$
such that $c'= \pi_N(X',Z')$ where $X'$ has distinct eigenvalues.
Fix such an $r$.
Since $\tau_{\beta(c')}(t)$ also has real coefficients, 
the above shows that $c' \in RC_n$. This means that 
$c' = \pi_n(X'', Z'')$ for some matrices $X'', Z''$ with real 
entries. Wilson's result \cite[Lemma 4.1]{W2} that $\beta$ 
intertwines the Calogero--Moser flows and the KP flows --
see \reref{flows} below -- implies that 
$c = \pi_n(X'' + r I_n, Z'') \in RC_n$, where $I_n$ 
is the identity matrix of size $n \times n$. 
(This also easily follows from \eqref{C-M-tau}.) 
This completes the proof of the theorem.
\end{proof}

\bre{flows} Wilson's adelic Grassmannian 
is invariant under the flows of the KP hierarchy, 
an infinite dimensional hamiltonian system 
on Sato's Grassmannian. On the level of tau functions 
the flows of the KP hierarchy account for shifts of the
variables $t_1, t_2, \ldots$ Similarly the Calogero--Moser 
spaces are the phase spaces of the Calogero--Moser 
hamiltonian systems, \cite{W2}. The latter are hamiltonian
systems with respect to a natural symplectic form on
$C_n$ and hamiltonians $h_k(\pi_n(X,Z))= (-1)^{k-1} \tr Z^k$.
They are the projections to $C_n$ of the linear flows
$(X,Z) \mapsto (X+kt(-Z)^{k-1}, Z)$ on $\ol{C}_n$. 
Wilson proved \cite[Lemma 4.1]{W2} that the map \eqref{C-M-tau}
intertwines the KP and the Calogero-Moser flows. 
Both the real loci of $C_n$ and the real locus of 
$\Gr^{\ad}$ are invariant under the Calogero-Moser flows and 
the KP flows, respectively, for real times.
\ere

We finish this section with a criterion for a point of Wilson's
adelic Grassmannian to belong to its real locus. Wilson's
Grassmannian possesses a remarkable involution called the bispectral
involution. It is defined by 
\[
\Psi_{b W} (x, z) = \Psi_W(z,x), \quad W \in \Gr^{ad}.
\]
This means that in \eqref{Psi} we expand all coefficients $a_j(x)$ 
for large $x$ and sum up the powers of $z$ in 
front of equal powers of $x$. 
It is a nontrivial statement proved in \cite{W1} that $\Psi_{b W}(x,z)$ has the form
\eqref{Psi}.
On level of the map \eqref{C-M-tau}, the bispectral 
involution is given by
\[
b(\beta(\pi_n(X,Z)))= \beta(\pi_n(Z^t, X^t)), 
\quad (X, Z) \in \ol{C}_n,
\]
see \cite{W2}.

\bth{Wad} Assume that $W \in \Gr^{\ad}$. If both 
$\tau_W(x, 0, \ldots)$ and $\tau_{bW}(x, 0, \ldots)$ 
have only real roots, then $W$ belongs to the real locus 
of $\Gr^{\ad}$.
\ethe

\thref{Wad} follows from \thref{2} and \thref{App}, 
and the fact that for $W = \beta(\pi_n(X,Z))$
\[
\tau_W(x, 0, \ldots) = (-1)^n \chi_X(-x), \quad
\tau_{bW}(x, 0, \ldots) = (-1)^n \chi_Z(-x),
\] 
where $\chi_X(x)$ denotes the characteristic polynomial
of $X$.
%%%%%%%%%%%%%%%%%%%%%%%%%%%%%%%%%%%%%%%%%%%%%%%%%%%%%%%%%%%%%%%
\sectionnew{Spaces of quasipolynomials} 
\label{quasi-pol}
Fix $\mu=(\mu_1,\dots,\mu_k) \in \Cset^k_{\reg}$, 
For a $k$-tuple of
finite dimensional subspaces $V_1, V_2, \dots, V_k$ of $\Cset[x]$,
we define the normalized Wronskian 
\[
\Wr(e^{\mu_1 x}V_1, \dots, e^{\mu_k x}V_k)
\]
as follows: choose a basis 
$\{q_1(x), \ldots, q_N(x)\}$ of 
$e^{\mu_1 x}V_1 \oplus \cdots \oplus e^{\mu_k x}V_k$ 
and set
\[
\Wr(e^{\mu_1 x}V_1, \dots, e^{\mu_k x}V_k)= 
a_0^{-1} e^{-(\mu_1 + \cdots + \mu_k) x} \Wr(q_1(x), \ldots, q_N(x))
\]
where $a_0$ is the leading coefficient of
\[
e^{-(\mu_1 + \cdots + \mu_k) x} \Wr(q_1(x), \ldots, q_N(x)).
\]
It is easy to check that the latter is a polynomial in $x$.
It is also straightforward to see that this definition does not depend on 
the choice of the basis $\{q_1(x), \ldots, q_N(x)\}$.

Recall the definition of Wilson's map $\gamma \colon G \to \Gr^{\ad}$
from \eqref{gamma}.

\ble{tau0} In the above setting
\[
\tau_{\gamma(\mu_1, \ldots, \mu_k, V_1, \ldots, V_k)}(x, 0, \ldots)
= s \Wr(e^{\mu_1 x}V_1, \dots, e^{\mu_k x}V_k), \quad 
s \in \Cset^*.
\]
\ele

This lemma can be extracted from \cite[Proposition 3.3]{SW}.
It follows by observing that the coefficient $a_1(x)$ of 
$\Psi_{\gamma(\mu_1, \ldots, \mu_k, V_1, \ldots, V_k)}(x,z)$ in
\eqref{Psi} is given by 
\[
a_1(x) = - \dx \log \tau_{\gamma(\mu_1, \ldots, \mu_k, V_1, \ldots, V_k)}
(x, 0, \ldots)
\]
because of \eqref{Sato},
and by 
\[
a_1(x) = - \dx \log \Wr(e^{\mu_1 x}V_1, \dots, e^{\mu_k x}V_k)
\]
because of \eqref{gamma}.

Using Wilson's map $\gamma$, we rederive the following theorem 
of Mukhin, Tarasov, and Varchenko \cite{MTV2} from \thref{2}. 

\bth{3} Let $(V_1, V_2, \dots, V_k)$ be a $k$-tuple of subspaces
of $\Cset[x]$. If $\mu_1, \ldots, \mu_k$ are real distinct numbers
and the polynomial $\Wr(e^{\mu_1 x}V_1, \dots, e^{\mu_k x}V_k)$ 
has only real roots, then each of the subspaces $V_1, \ldots, V_k$ has a 
basis consisting of polynomials with real coefficients.
\ethe

\begin{proof} Fix a collection 
such that 
$\Wr(e^{\mu_1 x}V_1, \dots, e^{\mu_k x}V_k)$ has only real roots and 
$\mu_1, \ldots, \mu_k$ are distinct real numbers.
We consider 
$(\mu_1, \ldots, \mu_k, V_1, \ldots, V_k)$
as an element of $G$. Let
\[
\gamma(\mu_1, \ldots, \mu_k, V_1, \ldots, V_k)= \beta(\pi_n(X,Z))
\]
for some $(X, Z) \in \ol{C}_n$. From \eqref{wave_beta} and 
\eqref{gamma}, we obtain that all eigenvalues of $Z$ are among 
$-\mu_1, \ldots, -\mu_k$, and thus are real. 
Comparing \eqref{C-M-tau} and \leref{tau0}, 
we see that 
\begin{equation}
\Wr(e^{\mu_1 x}V_1, \dots, e^{\mu_k x}V_k) = (-1)^n \chi_X(-x).
\label{tauWr}
\end{equation} 
Therefore all eigenvalues of $X$ are real 
as well. \thref{2} then implies that $\pi_n(X,Z) \in RC_n$ and 
\thref{App} implies that $\beta(\pi_n(X,Z))$ belongs to the real 
locus of $\Gr^{ad}$. Applying \thref{App} again, 
we obtain that 
\[
\gamma(\mu_1, \ldots, \mu_k, V_1, \ldots,  V_k)
=\gamma(\nu_1, \ldots, \nu_l, U_1, \ldots, U_l)
\]
for some $(\nu_1, \ldots, \nu_l, U_1, \ldots, U_k) \in G'$ 
with the properties of \thref{App}(5). Because of part 
(b) of \thref{Wilson}, $\{\nu_1, \ldots, \nu_l\}$ is a subset
of $\{\mu_1, \ldots, \mu_k\}$. So all $\nu_1, \ldots, \nu_l$ 
are real, and therefore each of the spaces 
$U_1, \ldots, U_l$ has a real basis because Condition 5
in \thref{App} is satisfied. \thref{Wilson} now implies that
the vector spaces $V_1, \ldots, V_k$ have the same properties.
This completes the proof of the theorem.
\end{proof} 

The Shapiro conjecture is the special 
case of \thref{3} when $k=1$ and $\mu_1=0$. 

Tracing back the relation between \thref{2}, 
\thref{Cher} and \thref{Wad}, 
we find three equivalent formulations of 
the Shapiro in the contexts of Calogero--Moser spaces, 
representations of $H_{0,1}(S_n)$ and Schur functions.

The tau functions of the KP hierarchy 
of the type $\tau_{\gamma(0,V)}(t)$ for $V \subset \Cset[x]$ are 
exactly the polynomial tau functions, i.e. those which are
polynomials and only depend on finitely many of the 
variables $t_1, t_2, \ldots$.
Recall that for a partition 
$\la=(\la_1 \geq \la_2 \geq \ldots \geq \la_l > 0)$ the corresponding 
Schur function is given by 
\begin{align}
\nn
&s_\la(p_1, \ldots, p_N) = \det (S_{\la_i + j-i}(p))_{i,j=1}^l, \\
\nn
&s_m(p_1, \ldots, p_m) = \sum_{j_1 + 2j_2 + \cdots=m}
\frac{p_1^{j_1}}{j_1} \frac{p_2^{j_2}}{j_2!} \dots,
\end{align}
where $N = \la_1 + l-1$. All polynomial tau functions of the KP 
hierarchy are of the type
\[
s_\la(t_1+ c_1, \ldots, t_N+ c_N)
\]
for some complex numbers $c_1, \dots, c_N \in \Cset$,
see \cite{SW,W2}. Finally, $\tau_{\beta(X,Z)}(t)$ is a polynomial 
tau function if and only if $Z$ is nilpotent, 
\cite[Proposition 6.1]{W2}. 

We obtain the following corollary for 
Schur functions.

\bco{Schur}
Fix a partition $\la=(\la_1 \geq \la_2 \geq \ldots \geq \la_l > 0)$
and an $N$-tuple of complex numbers $c_1, \ldots, c_N \in \Cset$, 
where $N = \la_1 + l-1$. If the corresponding Schur function
$s_\la(p_1, \ldots, p_N)$ has the property that
the polynomial $s_\la(x+ c_1, c_2, \ldots, c_N)$ has only real roots, 
then $c_1, c_2, \ldots, c_N \in \Rset$.
\eco

On the other hand we obtain that the Shapiro Conjecture is equivalent 
to the special case of \thref{Cher} for nilpotent actions of $y_1$, and 
to the following special case of \thref{2}:
\[
\Upsilon^{-1}( \Rset^{(n)} \times 0 ) \subset RC_n.
\] 
The latter equivalence was independently observed in \cite{MTV2}.
%%%%%%%%%%%%%%%%%%%%%%%%%%%%%%%%%%%%%%%%%%%%%%%%%%%%%%%%%%%%%%%%%%%%

%%%%%%%%%%%%%%%%%%%%%%%%%
\end{document}